\documentclass[a4paper]{article}

\usepackage[utf8]{inputenc}
\usepackage[T1]{fontenc}
\usepackage{lmodern}
\usepackage{amsmath,amsbsy,amssymb,bbm}
\usepackage{multirow,multicol,cite,enumerate,placeins}
\usepackage{graphicx}
\usepackage{color}
\usepackage[ruled]{algorithm2e}
\usepackage[hidelinks, bookmarks=false]{hyperref}

\usepackage[top=3cm,bottom=3cm,right=3.2cm,left=3.2cm]{geometry}

\newtheorem{lemma}{\bf Lemma}[section]
\newtheorem{theorem}[lemma]{\bf Theorem}
\newtheorem{cor}[lemma]{\bf Corollary}
\newtheorem{assumpt}[lemma]{\bf Assumption}

\def\PP{ {\rm I} \kern-.15em {\rm P} }

\def\Nb{\mbox{$\mathbb N$}}
\def\Rb{\mbox{$\mathbb R$}}
\def\Xb{\mbox{$\mathbb X$}}

\def\Bc{\mbox{$\mathcal B$}}

\def\Hc{\mbox{$\mathcal H$}}
\def\Nc{\mbox{$\mathcal N$}}
\def\Sc{\mbox{$\mathcal S$}}

\def \1{\mathbbm{1} }

\title{Improved bounds for Square-Root Lasso \\ and Square-Root Slope}
\author{Alexis Derumigny\thanks{ENSAE-CREST,
5, avenue Henry Le Chatelier,
TSA 96642, 91764 Palaiseau cedex, France. alexis.derumigny@ensae.fr}}
\date{\today}

\begin{document}

\maketitle

\begin{abstract}
    {\color{black} Extending the results of Bellec, Lecué and Tsybakov \cite{bellec2016slope} to the setting of sparse high-dimensional linear regression with unknown variance, we show that two estimators, the Square-Root Lasso and the Square-Root Slope can achieve the optimal minimax prediction rate, which is
    $(s/n) \log \left( p/s \right)$, up to some constant, under some mild conditions on the design matrix.}
    Here, $n$ is the sample size, $p$ is the dimension and $s$ is the sparsity parameter.
    We also prove optimality for the estimation error in the $l_q$-norm, with $q \in [1,2]$ for the Square-Root Lasso, and in the $l_2$ and sorted $l_1$ norms for the Square-Root Slope.
    Both estimators are adaptive to the unknown variance of the noise. The Square-Root Slope is also adaptive to the sparsity $s$ of the true parameter.
    {\color{black} Next, we prove that any estimator depending on $s$ which attains the minimax rate admits an adaptive to $s$ version still attaining the same rate.
    We apply this result to the Square-root Lasso.
    Moreover, for both estimators, we obtain valid rates for a wide range of confidence levels, and} improved concentration properties as in \cite{bellec2016slope} where the case of known variance is treated.
    Our results are non-asymptotic.
\end{abstract}

\medskip

\noindent
{\bf MCS:} Primary 62G08; secondary 62C20, 62G05.

\medskip


\bigskip

\section{Introduction}

In a recent paper by Bellec, Lecué and Tsybakov \cite{bellec2016slope}, it is shown that there exist high-dimensional statistical methods realizable in polynomial time that achieve the minimax optimal rate $(s/n) \log \left( p/s \right)$ in the context of sparse linear regression.
Here, $n$ is the sample size, $p$ is the dimension and $s$ is the sparsity parameter.
The result is achieved by the Lasso and Slope estimators, and the Slope estimator is adaptive to the unknown sparsity $s$.
Bounds for more general estimators are proved by Bellec, Lecué and Tsybakov \cite{bellec2016bounds, bellec2017towards}. These articles also establish bounds in deviation that hold for any confidence level and for the risk in expectation.
However, the estimators considered in \cite{bellec2016slope, bellec2016bounds, bellec2017towards} require the knowledge of the noise variance $\sigma^2$. To our knowledge, no polynomial-time methods, which would be at the same time optimal in a minimax sense and adaptive both to $\sigma$ and $s$ are available in the literature.

\medskip

Estimators similar to the Lasso, but adaptive to $\sigma$ are the Square-Root Lasso and the related Scaled Lasso, introduced by Sun and Zhang \cite{sun2012scaled} and Belloni, Chernozhukov and Wang \cite{belloni2011square}. It has been shown to achieve the rate $(s/n) \log(p)$ in deviation with the value of the tuning parameter depending on the confidence level.
A variant of this estimator is the Heteroscedastic Square-Root Lasso, which is studied in more general nonparametric and semiparametric setups by Belloni, Chernozhukov and Wang \cite{belloni2014pivotal}, but it also achieves the rate $(s/n) \log(p)$ and depends on the confidence level.
We refer to the book by Giraud \cite{giraud2014introduction} for the link between the Lasso and the Square-Root Lasso and a short proof of oracle inequalities for the Square-root Lasso.
In summary, there are two points to improve for the Square-root Lasso method:
\begin{enumerate}[(i)]
    \item The available results on oracle inequalities are valid only for the estimators depending on the confidence level.
    Thus, one cannot have an oracle inequality for one given estimator at any confidence level except the one that was used to design it.
    \item The obtained rate is $(s/n) \log(p)$ which is greater than the minimax rate $(s/n) \log(p/s)$.
\end{enumerate}

\medskip

The Slope, which is an acronym for Sorted L-One Penalized Estimation, is an estimator introduced by Bogdan et al. \cite{bogdan2015slope}, that is close to the Lasso, but uses the sorted $l_1$ norm instead of the standard $l_1$ norm for penalization.
Su and Candès \cite{su2016slope} proved that, as opposed to the Lasso, the Slope estimator is asymptotically minimax, in the sense that it attains the rate $(s/n) \log(p/s)$ for two isotropic designs, that is either for $\Xb$ deterministic with $\frac{1}{n}\Xb^T \Xb = I_{p \times p}$ or when $\Xb$ is a matrix with i.i.d. standard normal entries.
{\color{black} Moreover, their result has not only the optimal minimax rate, but also the exact optimal constant.}
General isotropic random designs are explored by Lecué and Mendelson \cite{lecue2016regularization}. For non-isotropic random designs and deterministic designs under conditions close to the Restricted Eigenvalue, the behavior of the Slope estimator is studied in \cite{bellec2016slope}.
The Slope estimator is adaptive only to $s$, and requires knowledge of $\sigma$, which is not available in practice.
In order to have an estimator which is adaptive both to $s$ and $\sigma$, we will use the Square-Root Slope, introduced by Stucky and van de Geer \cite{stucky2016sharp}. They give oracle inequalities for a large group of square-root estimators, including the new Square-Root Slope, but still following the scheme where (i) and (ii) cannot be avoided.
{\color{black} The square-root estimators are also members of a more general family of penalized estimators defined by Owen \cite[equations (8)-(9)]{owen2007robust} ; using their notation, these estimators correspond to the case where $\Hc_M$ is the squared loss and $\Bc_M$ is a norm (either the $l_1$ norm or the slope norm).}
\medskip

{\color{black} The paper is organized as follows. In Section \ref{section:framework}, we provide the main definitions and notations.} In Section \ref{section:optimal_sqrtLasso}, we show that the Square-Root Lasso is minimax optimal if $s$ is known while being adaptive to $\sigma$ under a mild condition on the design matrix (SRE). 
{\color{black} In Section \ref{section:adaptation_Lepski}, we show that any sequence of estimators can be made adaptive to the sparsity parameter $s$, while keeping the same rate up to some constant, with a computational cost increased by a factor of $\log(s_*)$ where $s_*$ is an upper bound on the sparsity parameter $s$. As an application, the Square-root Lasso modified by this procedure is still optimal while being now adaptive to $s$ (in addition of being already adaptive to $\sigma$).
In Section \ref{section:algorithms_sqrtSlope}, we show how to adapt any algorithm for computing the Slope estimator to the case of the Square-root Slope estimator.}
In Section \ref{section:optimal_sqrtSlope}, we study the Square-Root Slope estimator, and show that it is minimax optimal and adaptive both to $s$ and $\sigma$, under a slightly stronger condition (WRE). The (SRE) and (WRE) conditions have already been studied by Bellec, Lecué and Tsybakov \cite{bellec2016slope} and hold with high probability for a large class of random matrices.
Moreover, {\color{black}the inequalities we obtain for each estimator are valid for a wide range of confidence levels.}
Proofs are given in Section \ref{section:proofs}.

\medskip

\section{The framework}
\label{section:framework}

We use the notation $|\cdot|_q$ for the $l_q$ norm, with $1 \leq q \leq \infty$,
and $|\cdot|_0$ for the number of non-zero coordinates of a given vector.
For any $v \in \Rb^p$, and any set of coordinates $J$,
we denote by $v_J$ the vector $(v_j \1 \{i \in J\})_{i=1, \dots, p}$, where $\1$ is the indicator function.
We also define the empirical norm of a vector
$u = (u_1, \dots, u_n)$ as
$||u||_n^2 := \frac{1}{n} \sum_{i=1}^n u_i^2$.
For a vector $v \in \Rb^p$, we denote by $v_{(j)}$ the $j$-th largest component of $v$.
As a particular case, $|v|_{(j)}$ is the $j$-th largest component of the vector $|v|$ whose components are the absolute values of the components of $v$.
We use the notation $\langle \cdot , \cdot \rangle$ for the inner product with respect to the Euclidean norm and $(e_j)_{j=1, \dots,p}$ for the canonical basis in $\Rb^p$.

\medskip

Let $Y \in \Rb^n$ be the vector of observations
and let $\Xb \in \Rb^{n \times p}$ be the design matrix.
We assume that the true model is the following
\begin{equation}
    Y = \Xb \beta^* + \varepsilon.
    \label{model:linear_regression}
\end{equation}
Here $\beta^* \in \Rb^p$ is the unknown true parameter. We assume that $\varepsilon$ is the random noise, with values in $\Rb^n$, distributed as $\Nc(0,\sigma^2 I_{n \times n})$, where $I_{n \times n}$ is the identity matrix.
{\color{black} We denote by $\PP_{\beta^*}$ the probability distribution of $Y$ satisfying (\ref{model:linear_regression}).
In what follows, we define the set
$B_0(s) := \{ \beta^* \in \Rb^p: |\beta^*|_0 \leq s \}$.}
In the high-dimensional framework, we have typically in mind the case where $s$ is small,
$p$ is large and possibly $p \gg n$.

\medskip

We define two square-root type estimators of $\beta^*$: the Square-Root Lasso $\hat \beta^{SQL}$ and the Square-Root Slope $\hat \beta^{SQS}$ by the following relations
\begin{align}
    \label{definition_hat_beta_SQL}
    &\hat \beta^{SQL} \in \arg \min_{\beta \in \Rb^p}
    \left( \frac{1}{\sqrt{n}} |Y-\Xb \beta|_2 + \lambda |\beta|_1 \right),
    \\ &
    \label{definition_hat_beta_SQS}
    \hat \beta^{SQS} \in \arg \min_{\beta \in \Rb^p}
    \left( \frac{1}{\sqrt{n}} |Y-\Xb \beta|_2 + |\beta|_* \right),
\end{align}
where $\lambda > 0$ is a tuning parameter to be chosen,
and the sorted $l_1$ norm, $|\cdot|_*$, is defined for all $u \in \Rb^p$ by
$|u|_* = \sum_{i=1}^p \lambda_j |u|_{(j)},$
with tuning parameters $\lambda_1 \geq \dots \geq \lambda_p > 0$.

\section{Optimal rates for the Square-Root Lasso}
\label{section:optimal_sqrtLasso}

In this section, we derive oracle inequalities with optimal rate for the Square-Root Lasso estimator.
We will use the \emph{Strong Restricted Eigenvalue} (SRE) condition, introduced in \cite{bellec2016slope}.
For $c_0 > 0$ and $s\in \{ 1, \dots, p \}$, it is defined as follows,

\medskip

\noindent
$SRE(s,c_0)$ \textbf{condition :}
\emph{The design matrix $\Xb$ satisfies
$\max_{j=1, \dots,p} || \Xb e_j ||_n \leq 1$
and
\begin{equation}
    \kappa(s) := \min_{\delta \in C_{SRE}(s,c_0): \delta \neq 0}
    \dfrac{||\Xb \delta||_n}{|\delta|_2} > 0,
    \label{eq:def_kappa_SRE}
\end{equation}
where $C_{SRE}(s,c_0) := \{ \delta \in \Rb^p : |\delta|_1 \leq (1+c_0) \sqrt{s} |\delta|_2 \}$ is a cone in $\Rb^p$.
}

\bigskip

The condition $\max_{j=1, \dots,p} || \Xb e_j ||_n \leq 1$ is standard and corresponds to a normalization. It is shown in \cite[Proposition 8.1]{bellec2016slope} that the SRE condition is equivalent to the Restricted Eigenvalue (RE) condition of \cite{bickel2009simultaneous} if that is considered in conjunction with such a normalization.
{\color{black} By the same proposition, the RE condition is also equivalent to the $s$-sparse eigenvalue condition, which is satisfied with high probability for a large class of random matrices. It is the case, if for instance, $n \geq C s \log(ep/s)$ and the rows of $\Xb$ satisfies the small ball condition, which is very mild, see, e.g. \cite{bellec2016slope}.
}

{\color{black} Note that the minimum in (\ref{eq:def_kappa_SRE}) is the same as the minimum of the function $\delta \mapsto ||\Xb \delta||_n$ on the set $C_{SRE}(s,c_0) \cap \{ \delta \in \Rb^p : |\delta|_2 = 1\}$, which is a continuous function on a compact of $\Rb^p$, therefore this minimum is attained. When there is no ambiguity over the choice of $s$, we will just write $\kappa$ instead of $\kappa(s)$.}

\begin{theorem}
    \label{th:minimaxSqrtLasso}
    Let $s\in \{ 1, \dots, p \}$
    and assume that the $SRE(s,5/3)$ condition holds.
    Choose the following tuning parameter
    \begin{equation}
        \lambda = \gamma \sqrt{\frac{1}{n} \log \left( \frac{2p}{s} \right)},
        \label{cond:lambda_sqrtLasso}
    \end{equation}
    and assume that
    \begin{align}
        \gamma \geq 16 + 4\sqrt{2} \quad \text{and} \quad
        \frac{s}{n} \log \left(\frac{2p}{s}\right)
        \leq \frac{9 \kappa^2}{256 \gamma^2}.
        \label{cond:sqrt_Lasso_2}
    \end{align}
	{\color{black}
    Then, for every $\delta_0 \geq \exp(-n/4\gamma^2)$ and every $\beta^* \in \Rb^p$ such that $|\beta^*|_0 \leq s$,
    with $\PP_{\beta^*}$-probability at least $1 - \delta_0 - (1+e^2)e^{-n/24}$, we have
    \begin{align}
        ||\Xb (\hat \beta^{SQL} - \beta^*)||_n
        &\leq \sigma \max \left( \frac{C_1}{\kappa^2} 
        \sqrt{\frac{s}{n} \log \left(\frac{p}{s} \right)}
        \; , \; C_2 \sqrt{\frac{\log(1/\delta_0)}{n}}
        \right) ,
        \label{ineq_sql_n} \\
        |\hat \beta^{SQL} - \beta^*|_q
        &\leq \sigma \max \left( \frac{C_3}{\kappa^2} 
        s^{1/q} \sqrt{\frac{1}{n} \log \left( \frac{2p}{s} \right)}
        \; , \; C_4 s^{1/q-1} \sqrt{\frac{\log^2(1/\delta_0)}{n \log(2p/s)}}\right),
        \label{ineq_sql_q}
    \end{align}
    }
    where $1 \leq q \leq 2$, and $C_1>0, \, C_2>0, \, C_3>0, \, C_4>0$ are constants depending only on $\gamma$.
\end{theorem}

The values of the constants $C_1$, $C_2$, $C_3$ and $C_4$ in Theorem \ref{th:minimaxSqrtLasso} can be found in the proof, in Section \ref{proof:minimaxSqrtLasso}.
{\color{black}
Using the fact that $\kappa \leq 1$ and choosing $\delta_0 = (s/p)^s$, we get the following corollary of Theorem \ref{th:minimaxSqrtLasso}.
\begin{cor}
    Under the assumptions of Theorem \ref{th:minimaxSqrtLasso}, with $\PP_{\beta^*}$-probability at least $1 - (s/p)^s - (1+e^2)e^{-n/24}$, we have
    \begin{align*}
        ||\Xb (\hat \beta^{SQL} - \beta^*)||_n
        &\leq \frac{C_2}{\kappa^2} \sigma
        \sqrt{\frac{s}{n} \log \left(\frac{p}{s} \right)}, \\
        |\hat \beta^{SQL} - \beta^*|_q
        &\leq \frac{C_4}{\kappa^2} \sigma s^{1/q}
        \sqrt{\frac{1}{n} \log \left( \frac{2p}{s} \right)},
    \end{align*}
    where $1 \leq q \leq 2$.
    \label{cor:minimax_optimality_SQL}
\end{cor}
}

\medskip

{\color{black} Theorem \ref{th:minimaxSqrtLasso} and Corollary \ref{cor:minimax_optimality_SQL} give bounds that hold with high probability for both the prediction error and the estimation error in the $l_q$ norm, for every $q$ in $[1,2]$. Note that the bounds are best when the tuning parameter is chosen as small as possible, i.e. with $\gamma = 16 + 4\sqrt{2}$.
As shown in Section 7 of Bellec, Lecué and Tsybakov \cite{bellec2016slope}, the rates of estimation obtained in the latter corollary are optimal in a minimax sense on the set
$B_0(s) := \{ \beta^* \in \Rb^p: |\beta^*|_0 \leq s \}$.
We obtain the same rate of convergence as \cite{bellec2016slope} (see the paragraph after Corollary 4.3 in \cite{bellec2016slope}) up to some multiplicative constant.

The rate is also the same as in Su and Candès \cite{su2016slope}, but the framework is quite different: we obtain a non-asymptotic bound in probability whereas they consider asymptotic bounds in expectation (cf. Theorem 1.1 in \cite{su2016slope}) and in probability (Theorem 1.2) but without giving an explicit expression of the probability that their bound is valid.
Our result is non-asymptotic and valid when general enough conditions on $\Xb$ are satisfied whereas the result in \cite{su2016slope} is asymptotic as $n \to \infty$, and valid for two isotropic designs, that is either for $\Xb$ deterministic with $\frac{1}{n}\Xb^T \Xb = I_{p \times p}$ or when $\Xb$ is a matrix with i.i.d. standard normal entries.

Similarly to \cite{bellec2016slope}, for each tuning parameter $\gamma$, there is a wide range of levels of confidence $\delta_0$ under which the bounds of Theorem \ref{th:minimaxSqrtLasso} are valid. However, \cite{bellec2016slope} allows for an arbitrary small confidence level while in our case, there is a lower bound on the size of the confidence level under which the rate is obtained.
Note that this bound can be made arbitrary small by choosing a sample size $n$ large enough.

Note that the possible values chosen for the tuning parameter $\lambda$ are independent of the underlying standard deviation $\sigma$, which is unknown in practice.}
This gives an advantage for the Square-Root Lasso over other methods such as the ordinary Lasso.
Nevertheless, this estimator is not adaptive to the sparsity $s$, so that we need to know that $|\beta^*|_0 \leq s$ in order to be able to apply this result. {\color{black} In the following section, we suggest a procedure to make the Square-root Lasso adaptive to $s$ while keeping its optimality and adaptivity to $\sigma$.}

{\color{black}
\section{Adaptation to sparsity by a Lepski-type procedure}
\label{section:adaptation_Lepski}

Let $s_*$ be an integer in $\{2, \dots, p/e\}$.
We want to show that the Square-Root Lasso can also achieve the minimax optimal bound, adaptively to the sparsity $s$ on the interval $[1, s_*]$ (in addition of being already adaptive to $\sigma$).
Following \cite{bellec2016slope}, we will use aggregation of at most $\log_2(s_*)$ Square-Root Lasso estimators with different tuning parameters to construct an adaptive estimator $\tilde \beta$ of $\beta$ and at the same time an estimator $\tilde s$ of the sparsity $s$.

\medskip

In the following, we use the notation $\kappa_* := \kappa(2 s_*)$.
Note that $\kappa_* = \min_{s=1, \dots, 2 s_*} \kappa(s)$.
Indeed, the function $\kappa(\cdot)$ is decreasing, because the minimization (\ref{eq:def_kappa_SRE}) is done on spaces that are growing with $s$, in the sense of the inclusion.
We will assume that the condition $SRE(2 s_*, 5/3)$ holds and that 
$(2s_*/n) \log \big(2p/(2s_*) \big)
\leq 9 \kappa_*^2 / (256 \gamma^2)$.
The functions $b \mapsto (b/n) \log (2p/b)$ and $\kappa(\cdot)$ are respectively increasing (by Lemma \ref{lemma:choice_function_w}) and decreasing, so this ensures that the second part of condition (\ref{cond:sqrt_Lasso_2}) is satisfied for any $s=1, \dots, 2 s_*$.

\medskip

We can reformulate Corollary \ref{cor:minimax_optimality_SQL} as follows: for any $s=1, \dots, 2 s_*$ and any $\gamma \geq 16 + 4\sqrt{2}$
\begin{align}
    \sup_{\beta^* \in B_0(s)} \PP_{\beta^*}
    \left( ||\Xb (\hat \beta^{SQL}_{(s,\gamma)} - \beta^*)||_n
    \leq \frac{C_2(\gamma)}{\kappa_*^2} \sigma
    \sqrt{\frac{s}{n} \log \left(\frac{p}{s} \right)} \, \right) 
    \geq 1 - \left( \frac{s}{p} \right)^s - (1+e^2)e^{-n/24},
    \label{eq:SQL_sup_s}
\end{align}
denoting by $\hat \beta^{SQL}_{(s,\gamma)}$ the estimator (\ref{definition_hat_beta_SQL}) with the tuning parameter $\lambda_{(s,\gamma)}$ given by (\ref{cond:lambda_sqrtLasso}).
Replacing $s$ by $2s$ in equation (\ref{eq:SQL_sup_s}), we get that for any $s=1, \dots, s_*$ and any $\gamma \geq 16 + 4\sqrt{2}$,
\begin{align}
    \sup_{\beta^* \in B_0(2s)} \PP_{\beta^*}
    \left( ||\Xb (\hat \beta^{SQL}_{(2s,\gamma)} - \beta^*)||_n
    \leq \frac{C_2(\gamma)}{\kappa_*^2} \sigma
    \sqrt{\frac{2s}{n} \log \left(\frac{p}{2s} \right)} \, \right) 
    \geq 1 - \left( \frac{2s}{p} \right)^{2s} - (1+e^2)e^{-n/24}.
    \label{eq:SQL_sup_2s}
\end{align}
Remark that $\lambda_{(s,\gamma)}
= \gamma \sqrt{\frac{1}{n} \log \left( \frac{2p}{s} \right)}
= \tilde \gamma \sqrt{\frac{1}{n} \log \left( \frac{2p}{s} \right) - \frac{\log(2)}{n}}
= \lambda_{(2s,\tilde \gamma)}$ for some $\tilde \gamma > \gamma$.
As a consequence, $\hat \beta^{SQL}_{(s,\gamma)} = \hat \beta^{SQL}_{(2s,\tilde \gamma)}$ and we can apply Equation (\ref{eq:SQL_sup_2s}), replacing $\gamma$ by $\tilde \gamma$ and we get
\begin{align}
    \sup_{\beta^* \in B_0(2s)} \PP_{\beta^*}
    \left( ||\Xb (\hat \beta^{SQL}_{(s,\gamma)} - \beta^*)||_n
    \leq \frac{C_2(\tilde \gamma)}{\kappa_*^2} \sigma
    \sqrt{\frac{2s}{n} \log \left(\frac{p}{2s} \right)} \, \right) 
    \geq 1 - \left( \frac{2s}{p} \right)^{2s} - (1+e^2)e^{-n/24}.
    \label{eq:SQL_sup_2s_s}
\end{align}
Note that equations (\ref{eq:SQL_sup_s}) and (\ref{eq:SQL_sup_2s_s}) are the same as equations (5.2) and (5.4) in Bellec, Lecué and Tsybakov \cite{bellec2016slope}, taking $C_0 := \max \big(C_2(\gamma), C_2(\tilde \gamma) \big) / \kappa_*^2$, except that we have a supplementary term $- (1+e^2)e^{-n/24}$.
Similarly, we deduce from Corollary \ref{cor:minimax_optimality_SQL} that
\begin{align}
    \sup_{\beta^* \in B_0(s)} \PP_{\beta^*}
    \left( |\Xb (\hat \beta^{SQL}_{(s,\gamma)} - \beta^*)|_q
    \leq \frac{C_4(\gamma)}{\kappa_*^2} \sigma s^{1/q}
    \sqrt{\frac{s}{n} \log \left(\frac{2p}{s} \right)} \, \right) 
    \geq 1 - \left( \frac{s}{p} \right)^s - (1+e^2)e^{-n/24},
    \label{eq:SQL_sup_s_q} \\
    \sup_{\beta^* \in B_0(2s)} \PP_{\beta^*}
    \left( |\hat \beta^{SQL}_{(s,\gamma)} - \beta^*|_q
    \leq \frac{C_4(\tilde \gamma)}{\kappa_*^2} \sigma s^{1/q}
    \sqrt{\frac{2s}{n} \log \left(\frac{2p}{2s} \right)} \, \right) 
    \geq 1 - \left( \frac{2s}{p} \right)^{2s} - (1+e^2)e^{-n/24}.
    \label{eq:SQL_sup_2s_s_q}
\end{align}

\medskip

We describe now an algorithm to compute this adaptive estimator.
The idea is to use an estimator $\tilde s$ of $s$ which can be written as $\tilde s := 2^{\tilde m}$ for some positive data-dependent integer $\tilde m$.
We will use the notation $M := \max \{ m \in \Nb: 2^{m} \leq s_* \}$, so that the number of estimators we consider in the aggregation is $M$.

\medskip

The suggested procedure is detailed in Algorithm 1 below, with the distance $d(\beta,\beta') = ||\Xb(\beta - \beta')||_n$
or $d(\beta,\beta') = |\beta - \beta'|_q$ for $q \in [1,2]$.
It can be used for any family of estimators $(\hat \beta_{(s)})_{s=1, \dots s_*}$, and chooses the best one in terms of the distance $d(\cdot, \cdot)$, resulting in an aggregated estimator $\tilde \beta$.
Note that the weight function $w(\cdot)$ used in the algorithm cannot depend on $\sigma$ as in \cite{bellec2016slope}, i.e. to have the form $w(b) = C_0 \sigma \sqrt{(b/n) \log(p/b)}$ (respectively $w(b) = C_0 \sigma b^{1/q} \sqrt{(1/n) \log(p/b)}\,$), because we are looking for a procedure adaptive to $\sigma$. Therefore, we will remove $\sigma$ from $w$ and use an estimate $\hat \sigma$.

\medskip

\begin{algorithm}
\label{algo:adaptive_Lepski}
\SetAlgoLined
    \KwIn{a distance $d(\cdot,\cdot)$ on $\Rb^p$}
    \KwIn{a function $w(\cdot):[1, s_*] \rightarrow \Rb_+$ satisfying Assumption \ref{assumpt:function_w} }
    \KwIn{a family of estimators
    $\big( \hat \beta_{(s)} \big)_{s = 1, \dots, s_*}$}
    $M \leftarrow \lfloor \log_2(s_*) \rfloor$ \;
    \For{$m\leftarrow 1$ \KwTo $M+1$}{
        compute the estimator $\hat \beta_{(2^{m})}$ \;
    }
    compute $\hat \sigma \leftarrow ||Y - \Xb \hat \beta_{(2^{M+1})} ||_n$ \;
    compute the set $S_1 \leftarrow \left\{ m \in \{1, \dots, M \} :
    d \Big( \hat \beta_{(2^{k-1})} , \hat \beta_{(2^{k})} \Big)
    \leq 4 \hat \sigma C_0 w(2^{k}),
    \text{ for all } k \geq m \right\}$ \;
    \leIf{$S_1 \neq \emptyset$}{$\tilde m \leftarrow \min S_1$}
    {$\tilde m \leftarrow M$}
    \KwOut{$\tilde s \leftarrow 2^{\tilde m}$}
    \KwOut{$\tilde \beta \leftarrow \hat \beta_{(\tilde s)}$}
\caption{Algorithm for adaptivity.}
\end{algorithm}

\begin{assumpt}
    The function $w(\cdot): [1, s_*] \rightarrow \Rb_+$
    satisfies the following conditions:
    \begin{enumerate}
        \item $w(\cdot)$ is increasing on $[1, s_*]$ ;
        
        \item There exists a constant $C' > 0$ such that, for all $m=1, \dots, M$, we have $\sum_{k=1}^m w(2^k) \leq C' \cdot w(2^m)$ ;
        
        \item There exists a constant $C''> 0$ such that, for all $b=1, \dots, s_*$, $w(2b) \leq C'' w(b)$.
    \end{enumerate}
    \label{assumpt:function_w}
\end{assumpt}

\begin{assumpt}
    The family of estimators $(\hat \beta_{(s)})_{s=1, \dots, s_*}$ satisfies
    \begin{equation*}
        \sup_{\beta^* \in B_0(2s)}
        \PP_{\beta^*} \Big( \sigma/2 \leq \hat \sigma \leq \alpha \sigma \Big)
        \leq u_{n,p,M},
    \end{equation*}
    with a constant $\alpha > 0$, $\hat \sigma := ||Y - \Xb \hat \beta_{(2^{M+1})} ||_n$, and $u_{n,p,M} > 0$.
    \label{assumpt:bound_sigma_hat}
\end{assumpt}

\FloatBarrier

\begin{theorem}
    Let $s_*\in \{ 2, \dots, p/e \}$ and let $(\hat \beta_{(s)})_{s=1, \dots, s_*}$ be a collection of estimators satisfying Assumption \ref{assumpt:bound_sigma_hat} such that, for any $s = 1, \dots, s_*$,
    \begin{align}
        \sup_{\beta^* \in B_0(s)} \PP_{\beta^*}
        \left( d(\hat \beta_{(s)} \, , \, \beta^*)
        \leq C_0 \sigma w(s) \, \right) 
        \geq 1 - \left( \frac{s}{p} \right)^{s} - u_n,
        \label{eq:general_sup_s}
    \end{align}
    and
    \begin{align}
        \sup_{\beta^* \in B_0(2s)} \PP_{\beta^*}
        \left( d(\hat \beta_{(s)} \, , \, \beta^*)
        \leq C_0 \sigma w(2s) \, \right) 
        \geq 1 - \left( \frac{2s}{p} \right)^{2s} - u_n,
        \label{eq:general_sup_2s}
    \end{align}
    for a constant $C_0 > 0$, a function $w(\cdot): [1, s_*] \rightarrow \Rb_+$
    satisfying Assumption \ref{assumpt:function_w}, and $u_n > 0$.
    
    \medskip
    
    Then, there exists a constant $C_5$, depending on $C_0, C', C'', C_2, \kappa$ and $\alpha$ such that, for all $\beta^* \in B_0(s)$, the aggregated estimator $\tilde \beta$ satisfies:
    \begin{align*}
        &\PP_{\beta^*} \bigg( d(\tilde \beta, \beta^*)
        \leq C_5 \cdot \sigma w(s) \bigg)
        % \\ & \hspace{5cm}
        \geq 1 \, - \, 3 (\log_2(s_*) + 1)^2 \left( \left( \frac{2s}{p} \right)^{2s}
        + u_n \right) - u_{n,p,M}.
    \end{align*}
    Furthermore, 
    \begin{align*}
        &\PP_{\beta^*} \big( \tilde s \leq s \big) \geq 1 -
        2 (\log_2(s_*) + 1)^2 \left( \left( \frac{2s}{p} \right)^{2s}
        + u_n \right) - u_{n,p,M}.
    \end{align*}

    \label{th:proc_Lepski_adaptivity}
\end{theorem}

This theorem is proved in Section \ref{proof:th:proc_Lepski_adaptivity}.
In particular, it implies that when $\hat \beta_{(s)} = \hat \beta^{SQL}_{(s,\gamma)}$, the aggregated estimator $\tilde \beta$ has the same rate on $B_0(s)$ as the estimators with known $s$.
We detail it below.
The following lemmas proved in Sections \ref{proof:lemma:choice_function_w} and \ref{proof:lemma:validity_assumption_hat_sigma} assure that Theorem \ref{th:proc_Lepski_adaptivity} can be applied to the family $\hat \beta_{(s)} = \hat \beta^{SQL}_{(s,\gamma)}$.

\begin{lemma}
    Assumption \ref{assumpt:function_w} is satisfied with the choices $w(b) = \sqrt{(b/n) \log(p/b)}$ and 
    $w(b) = b^{1/q} \sqrt{(1/n) \log(2p/b)}$, for $q \in [1,2]$.
    \label{lemma:choice_function_w}
\end{lemma}

\begin{lemma}
    Assume that the $SRE(2s_*,5/3)$ condition holds and
    \begin{align*}
        \gamma \geq 16 + 4\sqrt{2} \quad \text{and} \quad
        \frac{2s_*}{n} \log \left(\frac{p}{s_*}\right)
        \leq 
        \min \left(\frac{9 \kappa_*^2}{256 \gamma^2} \, , \, \frac{\kappa_*^4}{2 C_2(\gamma)^2} \left( \frac{1}{\sqrt{2}} - \frac{1}{2} \right)^2 \right),
    \end{align*}
    where $\kappa_* := \kappa(2 s_*)$.
    Then Assumption \ref{assumpt:bound_sigma_hat} is satisfied with the choice $(\hat \beta_{(s)})_{s=1, \dots, s_*} = (\hat \beta^{SQL}_{(s,\gamma)})_{s=1, \dots, s_*}$, \\ $\alpha = 2 + \frac{3\sqrt{2} C_2(\gamma)}{16 \kappa \gamma}$ and $u_{n,p,M} = (2^{M+1}/p)^{2^{M+1}} - (1+e^2)e^{-n/24}$.
    \label{lemma:validity_assumption_hat_sigma}
\end{lemma}

\medskip

Combining equations (\ref{eq:SQL_sup_s}), (\ref{eq:SQL_sup_2s_s}) with Theorem \ref{th:proc_Lepski_adaptivity} and Lemmas \ref{lemma:choice_function_w} and \ref{lemma:validity_assumption_hat_sigma}, we obtain the following results for the case of the Square-root Lasso.

\begin{cor}
    Under the same assumptions as in Lemma \ref{lemma:validity_assumption_hat_sigma}, using Algorithm \ref{algo:adaptive_Lepski}, with $(\hat \beta_{(s)})_{s=1, \dots, s_*} = (\hat \beta^{SQL}_{(s,\gamma)})_{s=1, \dots, s_*}$, the distance $d(\beta,\beta') = ||\Xb(\beta - \beta')||_n$, and the weight $w(b) = \sqrt{(b/n) \log(p/b)}$, we have that, for all $\beta^* \in B_0(s)$, the aggregated estimator $\tilde \beta$ satisfies
    \begin{align*}
        &\PP_{\beta^*} \bigg( ||\Xb(\tilde \beta - \beta^*)||_n
        \leq C_5 \cdot \sigma \sqrt{\frac{s}{n} \log\left( \frac{p}{s} \right)} \, \bigg)
        % \\ & \hspace{5cm}
        \geq 1 \, - \, 3 (\log_2(s_*) + 1)^2 \left( \left( \frac{2s}{p} \right)^{2s}
        + u_n \right) - u_{n,p,M},
    \end{align*}
    and 
    \begin{align*}
        &\PP_{\beta^*} \big( \tilde s \leq s \big) \geq 1 -
        2 (\log_2(s_*) + 1)^2 \left( \left( \frac{2s}{p} \right)^{2s}
        + u_n \right) - u_{n,p,M},
    \end{align*}
    where $u_n = (1+e^2)e^{-n/24}$,
    $u_{n,p,M} = (2^{M+1}/p)^{2^{M+1}} - (1+e^2)e^{-n/24}$,
    and $C_5$ is a constant depending only on $\gamma$ and $\kappa_*$.
    \label{cor:sql_adaptive_s_norm_n}
\end{cor}

\medskip

\begin{cor}
    Under the same assumptions as in Lemma \ref{lemma:validity_assumption_hat_sigma}, using Algorithm \ref{algo:adaptive_Lepski}, with $(\hat \beta_{(s)})_{s=1, \dots, s_*} = (\hat \beta^{SQL}_{(s,\gamma)})_{s=1, \dots, s_*}$, the distance $d(\beta,\beta') = |\beta - \beta'|_q$, and the weight $w(b) = b^{1/q} \sqrt{(1/n) \log(2p/b)}$, for $q \in [1 ; 2]$, we have that, for all $\beta^* \in B_0(s)$, the aggregated estimator $\tilde \beta$ satisfies
    \begin{align*}
        &\PP_{\beta^*} \bigg( |\tilde \beta - \beta^*|_q
        \leq C_5 \cdot \sigma s^{1/q} \sqrt{\frac{1}{n} \log\left( \frac{p}{s} \right)} \, \bigg)
        % \\ & \hspace{5cm}
        \geq 1 \, - \, 3 (\log_2(s_*) + 1)^2 \left( \left( \frac{2s}{p} \right)^{2s}
        + u_n \right) - u_{n,p,M},
    \end{align*}
    and 
    \begin{align*}
        &\PP_{\beta^*} \big( \tilde s \leq s \big) \geq 1 -
        2 (\log_2(s_*) + 1)^2 \left( \left( \frac{2s}{p} \right)^{2s}
        + u_n \right) - u_{n,p,M},
    \end{align*}
    where $u_n = (1+e^2)e^{-n/24}$,
    $u_{n,p,M} = (2^{M+1}/p)^{2^{M+1}} - (1+e^2)e^{-n/24}$,
    and $C_5$ is a constant depending only on $\gamma$ and $\kappa_*$.
    \label{cor:sql_adaptive_s_norm_q}
\end{cor}

Thus, we have shown that the suggested aggregated procedure based on the Square-root Lasso is adaptive to $s$ while still being adaptive to $\sigma$ and minimax optimal. Note that the computational cost is multiplied by $O(\log(s_*))$.

\bigskip

}

{\color{black}

\section{Algorithms for computing the Square-root Slope}
\label{section:algorithms_sqrtSlope}

In this part, our goal is to provide algorithms for computing the square-root Slope estimator. A natural idea is revisiting the algorithms used for the square-root Lasso and for the Slope, then adapting or combining them.

\medskip

Belloni, Chernozhukov and Wang \cite[Section 4]{belloni2011square} have proposed to compute the Square-root Lasso estimator by reducing its definition to an equivalent problem, which can be solved by interior-point or first-order methods. 
The equivalent formulation as the Scaled Lasso, introduced by Sun and Zhang \cite{sun2012scaled} allows one to view it as a joint minimization in $(\beta, \sigma)$. Sun and Zhang \cite{sun2012scaled} propose an iterative algorithm which alternates estimation of $\beta$ using the ordinary Lasso and estimation of $\sigma$.

\medskip

Zeng and Figueiredo \cite{zeng2014ordered} studied several algorithms related to estimation of the regression with the ordered weighted $l_1$-norm, which is the Slope penalization. Bogdan et al. \cite{bogdan2015slope} provide an algorithm for computing the Slope estimator using a proximal gradient.

\medskip

As in the case of the Square-root Lasso, we still have for any $\beta$,
\begin{equation}
    ||Y - \Xb \beta ||_n = \min_{\sigma > 0}
    \left( \sigma + \frac{||Y - \Xb \beta ||_n^2}{\sigma} \right),
\end{equation}
where the minimum is attained for $\hat \sigma = ||Y - \Xb \beta ||_n$.
As a consequence,
\begin{align*}
    \hat \beta^{SQS} \in \arg \min_{\beta \in \Rb^p}
    \big(||Y-\Xb \beta||_n + |\beta|_* \big)
\end{align*}
is equivalent to take the estimator $\hat \beta$ in the joint minimization program
\begin{align*}
    (\hat \beta, \hat \sigma) \in
    \underset{\beta \in \Rb^p, \, \sigma > 0}{\arg \min}
    \left(\sigma + \frac{||Y - \Xb \beta ||_n^2}{\sigma} + |\beta|_* \right).
\end{align*}

Alternating minimization in $\beta$ and in $\sigma$ gives an iterative procedure for a "Scaled Slope" (see Algorithm \ref{algo:scaled_Slope}).

\begin{algorithm}
\label{algo:scaled_Slope}
\SetAlgoLined
    \KwIn{explained variable $Y$, design matrix $\Xb$ ;}
    \KwIn{tuning parameters $\lambda_1 \leq \dots \leq \lambda_p$ ;}
    choose some initialization value for $\hat \sigma$, for example the standard deviation of $Y$ \;
    \Repeat(){convergence }{
    estimate $\hat \beta$ by the Slope algorithm with the parameters
    $\hat \sigma \cdot \lambda_1, \dots, \hat \sigma \cdot \lambda_p$ \;
    estimate $\hat \sigma$ by $||Y - \Xb \hat \beta ||_n$ ;
    }
    \KwOut{a joint estimator $\big( \hat \beta , \, \hat \sigma \big)$ ;}
    \caption{Scaled Slope algorithm}
\end{algorithm}

}

\section{Optimal rates for the Square-Root SLOPE}
\label{section:optimal_sqrtSlope}

In this part, we will use another condition, the \emph{Weighted Restricted Eigenvalue} condition, introduced in \cite{bellec2016slope}.
For $c_0 > 0$ and $s\in \{ 1, \dots, p \}$, it is defined as follows,

\medskip

\noindent
$WRE(s,c_0)$ \textbf{condition :}
\emph{The design matrix $\Xb$ satisfies
$\max_{j=1, \dots,p} || \Xb e_j ||_n \leq 1$
and
\begin{equation}
    \kappa' := \min_{\delta \in C_{WRE}(s,c_0): \delta \neq 0}
    \dfrac{||\Xb \delta||_n}{|\delta|_2} > 0,
\end{equation}
where $C_{WRE}(s,c_0) := \{ \delta \in \Rb^p : |\delta|_* \leq (1+c_0) |\delta|_2 \sqrt{\sum_{j=1}^s \lambda_j^2} \}$ is a cone in $\Rb^p$.
}

\bigskip

To obtain the following result, we assume that the Weighted Restricted Eigenvalue condition holds. This condition is shown to be only slightly more constraining than the usual Restricted Eigenvalue condition of \cite{bickel2009simultaneous}, but is nevertheless satisfied with high probability for a large class of random matrices, see Bellec, Lecué and Tsybakov \cite{bellec2016slope} for a discussion.
{\color{black} Note that, in a similar way as in definition (\ref{eq:def_kappa_SRE}), the minimum is attained. Indeed, $\kappa'$ is equal to the minimum of the function $\delta \mapsto ||\Xb \delta||_n$ on the set $C_{WRE}(s,c_0) \cap \{ \delta \in \Rb^p : |\delta|_2 = 1\}$, which is a continuous function on a compact of $\Rb^p$.}

\begin{theorem}
    \label{th:minimaxSqrtSlope}
    Let $s\in \{ 1, \dots, p \}$
    and assume that the $WRE(s,20)$ condition holds.
    Choose the following tuning parameters
    \begin{equation}
        \lambda_j = \gamma' \sqrt{ \frac{\log(2p/j)}{n} },
        \text{ for } j=1, \dots, p,
        \label{cond:lambda_sqrtSlope}
    \end{equation}
    and assume that
    \begin{align}
        \gamma' \geq 16 + 4\sqrt{2} \quad \text{and} \quad
        \frac{s}{n} \log \left(\frac{2ep}{s} \right)
        \leq \frac{\kappa'{}^2}{256 \gamma'{}^2}.
        \label{cond:sqrt_Slope_2}
    \end{align}
    Then, for every $\delta_0 \geq \exp(-n/4\gamma'{}^2)$ and every $\beta^* \in \Rb^p$ such that $|\beta^*|_0 \leq s$,
    with $\PP_{\beta^*}$-probability at least $1 - \delta_0 - (1+e^2)e^{-n/24}$, we have
    {\color{black}
    \begin{align}
        ||\Xb (\hat \beta^{SQS} - \beta^*)||_n
        & \leq \sigma \max \Bigg(
        \frac{C'_1}{\kappa'} \sqrt{\frac{s}{n} \log \left(\frac{p}{s} \right)}
        \, , \, C'_2 \sqrt{\frac{\log(1/\delta_0)}{n}} \Bigg),
        \\
        |\hat \beta^{SQS} - \beta^*|_*
        & \leq \sigma \max \Bigg(
        \frac{C'_1}{\kappa'{}^2} \frac{s}{n} \log \left(\frac{p}{s} \right)
        \, , \, C'_2 \frac{\log(1/\delta_0)}{n} \Bigg),
        \\
        |\hat \beta^{SQS} - \beta^*|_2
        & \leq \sigma \max \Bigg(
        \frac{C'_1}{\kappa'{}^2} \sqrt{\frac{s}{n} \log \left(\frac{p}{s} \right)}
        \, , \, C'_2 \sqrt{\frac{\log^2(1/\delta_0)}{sn\log(p/s)}} \Bigg),
    \end{align}
    }
    for constants $C'_1>0$ and $C'_2>0$ depending only on $\gamma'$.
\end{theorem}

\medskip

The values of the constants $C'_1$ and $C'_2$ can be found in the proof, in Subsection \ref{proof:minimaxSqrtSlope}.
{\color{black} Note that the bounds are best when the tuning parameters is chosen as small as possible, i.e. using the choice $\gamma' = 16 + 4\sqrt{2}$.}
{\color{black}
Using the fact that $\kappa' \leq 1$ and choosing $\delta_0 = (s/p)^s$, we get the following corollary.
\begin{cor}
    Under the assumptions of Theorem \ref{th:minimaxSqrtSlope}, with $\PP_{\beta^*}$-probability at least $1 - (s/p)^s - (1+e^2)e^{-n/24}$, we have
    \begin{align*}
        ||\Xb (\hat \beta^{SQS} - \beta^*)||_n
        & \leq \frac{C'_1}{\kappa'} \sigma
        \sqrt{\frac{s}{n} \log \left(\frac{p}{s} \right)},
        \\
        |\hat \beta^{SQS} - \beta^*|_*
        & \leq \frac{C'_1}{\kappa'{}^2} \sigma
        \frac{s}{n} \log \left(\frac{p}{s} \right),
        \\
        |\hat \beta^{SQS} - \beta^*|_2
        & \leq \frac{C'_1}{\kappa'{}^2} \sigma
        \sqrt{\frac{s}{n} \log \left(\frac{p}{s} \right)},
    \end{align*}
\end{cor}
}
These results show that the Square-Root Slope estimator, with a given choice of parameters, attains the optimal rate of convergence in the prediction norm $||\cdot||_n$ and in the estimation norm $|\cdot|_2$.
We also provide a bound on the sorted $l_1$ norm $|\cdot|_*$ of the estimation error.
One can note that the choice of $\lambda_i$ that allows us to obtain optimal bounds does not depend on the level of confidence $\delta_0$, {\color{black}but only influence the size of the range of valid $\delta_0$}.
This improves upon the oracle result of Stucky and van de Geer \cite{stucky2016sharp}, in which the parameter does depend on the level of confidence and the rate does not scale in the optimal way, i.e., as $\sqrt{(s/n) \log(p/s)}$.
Moreover, we can see that our estimator is independent of the underlying standard deviation $\sigma$ and of the sparsity $s$, even if the rates depend on them.
{\color{black} Note that, up to some multiplicative constant, we obtain the same rates as for the Slope in Bellec, Lecué and Tsybakov \cite{bellec2016slope}. In Su and Candès \cite{su2016slope}, the Slope estimator is proved to attain the sharp constant in the asymptotic framework where $\sigma$ is known and for specific $\Xb$ ; whereas here we obtain only the minimax rates, but in a non-asymptotic framework, and under general assumptions on the design matrix $\Xb$.
}

\medskip

For this estimator, we did not provide a bound for the $l_1$ norm, for the same reasons as in \cite{bellec2016slope}. Indeed, the coefficients $\lambda_j$ of the components of $\beta$ are different in the sorted norm.
As a consequence, we do not provide inequalities for $l_q$ norms when $q<2$, that are obtained by interpolation between the $l_1$ and $l_2$ norms.

\section{Proofs}
\label{section:proofs}

\subsection{Preliminary lemmas}

Let $\beta^* \in \Rb^p$, $\Sc \subset \{1,\dots,p\}$ with cardinality $s$ and denote by $\Sc^C$ the complement of $\Sc$. For $i \in \{1,\dots,p\}$, let $\beta^*_i$ be the $i$-th component of $\beta^*$ and assume that for every $i \in \Sc^C$, $\beta^*_i = 0$.

\medskip

%==========================

\begin{lemma} We have
    $|(\hat \beta^{SQL} - \beta^*)_{\Sc^C}|_1
    \leq |(\hat \beta^{SQL} - \beta^*)_{\Sc}|_1
    + \dfrac{1}{\lambda \sqrt{n} | \varepsilon |_2}
    \left \langle \Xb^T \varepsilon \; , \; \hat \beta^{SQL} - \beta^* \right \rangle.$
    \label{equ_beta_mc}
\end{lemma}
The proof follows from the arguments in Giraud \cite[pages 110-111]{giraud2014introduction}, and it is therefore omitted.

%==========================
\bigskip

\begin{lemma}
    Let $u \in \Rb^p$ be defined by $u := \hat \beta^{SQS} - \beta^*$.
    We have $$\sum_{j=s+1}^p \lambda_j |u|_{(j)}
    \leq \sum_{j=1}^s \lambda_j |u|_{(j)}
    + \frac{1}{\sqrt{n} | \varepsilon |_2}
    \left \langle \Xb^T \varepsilon \; , \; u \right \rangle.$$
    \label{equ_beta_sqs}
\end{lemma}
{\it Proof :} We combine the arguments from Giraud \cite[pages 110-111]{giraud2014introduction}, and from the proof of Lemma A.1 in \cite{bellec2016slope}.
First, we remark that the sorted $l_1$ norm can be written as follows, for any $v \in \Rb^p$,
\begin{equation*}
    |v|_* = \max_\phi \sum_{j=1}^p \lambda_j \left| v_{\phi(j)} \right|,
\end{equation*}
where the maximum is taken over all permutations $\phi = (\phi(1), \dots, \phi(p))$ of $\{1, \dots, p \}$.

\medskip

\noindent
By definition,
$\hat \beta^{SQS}$ is a minimizer of (\ref{definition_hat_beta_SQS}),
so we have
\begin{align*}
    |Y-\Xb \hat \beta^{SQS}|_2 - |Y-\Xb \beta^*|_2
    &\leq \sqrt{n} \left( |\beta^*|_* - | \hat \beta^{SQS}|_* \right).
\end{align*}
Let $\phi$ be any permutation of $\{1, \dots, p \}$ such that
\begin{equation}
    |\beta^*|_* = \sum_{j=1}^s \lambda_j |\beta^*_{\phi(j)}| \, \text{ and } \,
    |u_{\phi(s+1)}| \geq |u_{\phi(s+2)}| \geq \dots \geq |u_{\phi(p)}|.
    \label{cond:permutation_phi}
\end{equation}
We have 
\begin{align*}
    |\beta^*|_* - | \hat \beta^{SQS}|_*
    &\leq \sum_{j=1}^s \lambda_j \Big( \big|\beta^*_{\phi(j)}\big|
    - \big| \hat \beta_{\phi(j)}^{SQS} \big| \Big)
    - \sum_{j=s+1}^p \lambda_j \big| \hat \beta_{\phi(j)}^{SQS} \big| \\
    &\leq \sum_{j=1}^s \lambda_j \big| u_{\phi(j)} \big|
    - \sum_{j=s+1}^p \lambda_j \big| \hat \beta_{\phi(j)}^{SQS} \big|
    = \sum_{j=1}^s \lambda_j \big| u_{\phi(j)} \big|
    - \sum_{j=s+1}^p \lambda_j \big| u_{\phi(j)} \big|.
\end{align*}
Since the sequence $\lambda_j$ is non-increasing, we have
$\sum_{j=1}^s \lambda_j |u_{\phi(j)}| \leq \sum_{j=1}^s \lambda_j |u|_{(j)}$.
The permutation $\phi$ satisfies (\ref{cond:permutation_phi}), therefore,
$\sum_{j=s+1}^p \lambda_j |u|_{(j)} \leq \sum_{j=s+1}^p \lambda_j |u_{\phi(j)}|$.
From the previous inequalities, we get that
\begin{align}
    |Y-\Xb \hat \beta^{SQS}|_2 - |Y-\Xb \beta^*|_2
    &\leq \sqrt{n} \left( \sum_{j=1}^s \lambda_j |u|_{(j)} - \sum_{j=s+1}^p \lambda_j |u|_{(j)} \right).
    \label{borne_sup_diff_SQS}
\end{align}
By convexity of the mapping $\beta \mapsto ||Y-X \beta||_2$, we have
\begin{equation}
    |Y-\Xb \hat \beta^{SQS}|_2 - |Y-\Xb \beta^*|_2
    \geq - \left \langle
    \frac{\Xb^T \varepsilon}{| \varepsilon |_2} \; , \; \hat \beta^{SQS} - \beta^* \right \rangle
    = -\frac{1}{| \varepsilon |_2}
    \left \langle \Xb^T \varepsilon \; , \; \hat \beta^{SQS} - \beta^* \right \rangle.
    \label{borne_inf_diff_SQS}
\end{equation}
Combining (\ref{borne_sup_diff_SQS}) and (\ref{borne_inf_diff_SQS}), we get
\begin{equation*}
    -\frac{1}{|\varepsilon |_2}
    \left \langle \Xb^T \varepsilon \; , \; \hat \beta^{SQS} - \beta^* \right \rangle
    \leq
    \sqrt{n}
    \left(\sum_{j=1}^s \lambda_j |u|_{(j)} - \sum_{j=s+1}^p \lambda_j |u|_{(j)} \right),
\end{equation*}
which concludes the proof.

\begin{flushright}
    $\Box$
\end{flushright}

%==========================
\bigskip

\begin{lemma} We have
    $|\Xb (\hat \beta^{SQL} - \beta^*)|_2^2
    \leq
    \left \langle \Xb^T \varepsilon \; , \; \hat \beta^{SQL} - \beta^* \right \rangle
    + \lambda \sqrt{n} |Y - \Xb \hat \beta^{SQL}|_2
    |\hat \beta^{SQL} - \beta^*|_1.$
    \label{first_order_consequence_sql}
\end{lemma}
\begin{lemma} We have
    $|\Xb (\hat \beta^{SQS} - \beta^*)|_2^2
    \leq
    \left \langle \Xb^T \varepsilon \; , \; \hat \beta^{SQS} - \beta^* \right \rangle
    + \sqrt{n} |Y - \Xb \hat \beta|_2
    |\hat \beta^{SQS} - \beta^*|_*.$
    \label{first_order_consequence_sqs}
\end{lemma}

\noindent
{\it Proof :} We will give a general proof of Lemmas \ref{first_order_consequence_sql} and \ref{first_order_consequence_sqs} in the case of an estimator defined by
\begin{align}
    \hat \beta := \arg \min_{\beta \in \Rb^p}
    \left( \frac{1}{\sqrt{n}} |Y-\Xb \beta|_2 + ||\beta|| \right),
\end{align}
where $||\cdot||$ is a norm on $\Rb^p$.
Lemmas \ref{first_order_consequence_sql} and \ref{first_order_consequence_sqs} are obtained as special cases corresponding to $||\cdot||=\lambda|\cdot|_1$ and $||\cdot||=|\cdot|_*$.
Denote by $||\cdot||_{dual}$ the norm dual to $||\cdot||$.

\medskip

Since $\hat \beta$ is optimal, we know that
$\Xb^T (Y - \Xb \hat \beta) / (\sqrt{n} |Y - \Xb \hat \beta|_2)$
belongs to the subdifferential of the function $||\cdot||$ evaluated at $\hat \beta$.
Thus, there exists $v \in \Rb^p$ such that $||v||_{dual} \leq 1$ and
\begin{align*}
    \frac{\Xb^T (Y - \Xb \hat \beta)}
    {\sqrt{n} |Y - \Xb \hat \beta|_2} + v = 0.
\end{align*}
Thus, we have
\begin{gather*}
    % \Xb^T(Y-\Xb \hat \beta)
    % + \lambda \sqrt{n} v |Y - \Xb \hat \beta|_2 = 0,   \\
    %
    %\Xb^T(\Xb \beta + \varepsilon -\Xb \hat \beta) + \sqrt{n} \lambda v |Y - \Xb \hat \beta|_2 = 0,   \\
    %
    % \langle \Xb(\beta - \hat \beta) \; , \;
    % \Xb(\beta - \hat \beta) \rangle
    % + \langle \Xb^T \varepsilon \; , \;
    % \beta - \hat \beta \rangle
    % + \lambda \sqrt{n} |Y - \Xb \hat \beta|_2 
    % \langle v \; , \; \beta - \hat \beta \rangle = 0, \\
    %
    |\Xb (\hat \beta - \beta^*)|_2^2 =
    \left \langle \Xb^T \varepsilon \; , \; \hat \beta - \beta^* \right \rangle
    + \sqrt{n} |Y - \Xb \hat \beta^*|_2
    \langle v \; , \; \hat \beta - \beta^* \rangle.
\end{gather*}
The conclusion results from the inequality
\begin{equation*}
    \langle v \; , \; \hat \beta - \beta^* \rangle
    \leq ||v||_{dual} ||\hat \beta - \beta^*||
    \leq ||\hat \beta - \beta^*||.
\end{equation*}

\begin{flushright}
    $\Box$
\end{flushright}

%==========================
\bigskip

\begin{lemma}
   We have $\gamma' \sqrt{(s/n) \log(2p/s)}
   \leq \sqrt{\sum_{j=1}^s \lambda_j^2}
    \leq \gamma' \sqrt{(s/n) \log(2ep/s)}$.
   \label{lemma:bound_sum_lambda}
\end{lemma}
{\it Proof :}
From Stirling's formula, we deduce that
$s \log(s/e) \leq \log(s!) \leq s \log(s)$. Therefore
\begin{equation*}
    s\log(2p/s) \leq \sum_{j=1}^s \log(2p/j)
    = \log(2p) - \log(s!) \leq s \log(2ep/s).
\end{equation*}
The conclusion follows from the definition of the $\lambda_j$ in (\ref{cond:lambda_sqrtSlope}).

\begin{flushright}
    $\Box$
\end{flushright}

%==========================
\bigskip

The following simple property is proved in Giraud \cite[page 112]{giraud2014introduction}. For convenience, it is stated here as a lemma.
\begin{lemma}
    With $\PP_{\beta^*}$-probability at least $1 - (1+e^2)e^{-n/24}$, we have
    \begin{equation*}
        \frac{\sigma}{\sqrt{2}}
        \leq \frac{|\varepsilon|_2}{\sqrt{n}}
        \leq 2 \sigma.
    \end{equation*}
    \label{lemma_epsilon_sigma}
\end{lemma}
We will also use the following theorem from Bellec, Lecué and Tsybakov
\cite[Theorem 4.1]{bellec2016slope}.

\begin{lemma}
    Let $0 < \delta_0 < 1$ and let $\Xb$ in $\Rb^{n \times p}$ be a matrix such that
    $\max_{j=1, \dots, p} || \Xb e_j ||_n \leq 1$. For any $u = (u_1, \dots u_p)$ in $\Rb^p$, we define :
    \begin{equation*}
        G(u) := (4 + \sqrt{2}) \sigma
        \sqrt{\frac{\log(1/\delta_0)}{n}} ||\Xb u||_n,
        \quad  
        H(u) := (4 + \sqrt{2}) \sum_{j=1}^p |u|_{(j)} \sigma \sqrt{\frac{\log(2p/j)}{n}},
    \end{equation*}
    
    \begin{equation*}
        \text{and} \quad F(u) := (4 + \sqrt{2}) \sigma \sqrt{\frac{\log(2p/s)}{n}}
        \left( \sqrt{s}|u|_2 + \sum_{j=s+1}^p |u|_{(j)} \right).
    \end{equation*}
    If $\varepsilon \sim \Nc(0,\sigma^2  I_{n \times n})$, then the random event
    \begin{equation*}
        \left\{ \frac{1}{n} \varepsilon^T \Xb u \leq \max \Big( H(u), G(u) \Big),
        \forall u \in \Rb^p \right\},
    \end{equation*}
    is of probability at least $1 - \delta_0/2$.
    
    \medskip
    
    Moreover, by the Cauchy-Schwarz inequality, we have
    $H(u) \leq F(u)$, for all $u$ in $\Rb^p$.
    \label{lemma_Th_4_1_BLT}
\end{lemma}

\subsection{Proof of Theorem \ref{th:minimaxSqrtLasso} }
\label{proof:minimaxSqrtLasso}

Lemma \ref{lemma_Th_4_1_BLT} allows one to control the random variable $\varepsilon^T \Xb u$ that appears in Lemmas \ref{equ_beta_mc}
and \ref{first_order_consequence_sql} with $u := \hat \beta^{SQL} - \beta^*$.
Our calculations will take place on an event of probability at least $1 - \delta_0 - (1+e^2)e^{-n/24}$, where both Lemmas \ref{lemma_epsilon_sigma} and \ref{lemma_Th_4_1_BLT} can be used.
Applying Lemma \ref{lemma_Th_4_1_BLT}, we will distinguish between the two cases : $G(u) \leq F(u)$ and $F(u) < G(u)$.

\bigskip

\noindent
{\it First case :} $G(u) \leq F(u)$.

\medskip

\noindent
Then we have
\begin{equation*}
    (4 + \sqrt{2}) \sqrt{\frac{\log(1/\delta_0)}{n}} ||\Xb u||_n
    \leq (4 + \sqrt{2}) \sqrt{\frac{\log(2p/s)}{n}} \left( \sqrt{s}|u|_2 + \sum_{j=s+1}^p |u|_{(j)} \right).
\end{equation*}

\medskip

\noindent
We will show first that $u$ is in the SRE cone, so that we can use the SRE assumption.
From Lemma \ref{equ_beta_mc}, we have
\begingroup \allowdisplaybreaks
\begin{align*}
    |u_{\Sc^C}|_1
    &\leq |u_{\Sc}|_1
    + \frac{1}{\lambda \sqrt{n} |\varepsilon|_2}
    \left \langle \Xb^T \varepsilon \; , \; \hat \beta^{SQL} - \beta^* \right \rangle \\
    &\leq |u_{\Sc}|_1
    + \frac{1}{\sqrt{n} \lambda |\varepsilon|_2}
    n \sigma (4 + \sqrt{2}) \sqrt{\frac{\log(2p/s)}{n}} \left( \sqrt{s}|u|_2 + \sum_{j=s+1}^p |u|_{(j)} \right) \\
    &\leq |u_{\Sc}|_1
    + \frac{1}{4}
    \left( \sqrt{s}|u|_2 + |u_{\Sc^C}|_1 \right),
\end{align*}
where in the last inequality, we have used Lemma \ref{lemma_epsilon_sigma} and assumption (\ref{cond:sqrt_Lasso_2}). We deduce that
% \begin{align*}
%     \frac{3}{4} |u_{\Sc^C}|_1
%     &\leq |u_{\Sc}|_1 + \frac{1}{4} \sqrt{s}|u|_2.
% \end{align*}
% As a consequence,
\begin{align*}
    \frac{3}{4} |u|_1
    &\leq \frac{7}{4} |u_{\Sc}|_1
    + \frac{1}{4} \sqrt{s} |u|_2
    \leq \frac{7}{4} \sqrt{s} |u|_2
    + \frac{1}{4} \sqrt{s} |u|_2 = 2 \sqrt{s} |u|_2.
\end{align*}
\endgroup
Therefore, we have
\begin{align}
    \label{condition_for_SRE}
    |u|_1 \leq \frac{8}{3} \sqrt{s} |u|_2,
\end{align}
and thus, the following inequality holds
$|u|_1 \leq (1+c_0) \sqrt{s} |u|_2$,
with $c_0 = 5/3$, allowing us to use the $SRE(s,5/3)$ assumption.

\medskip

\noindent
From Lemmas \ref{first_order_consequence_sql} and \ref{lemma_Th_4_1_BLT}, and using that, in view of the $SRE(s,5/3)$ condition, $||\Xb u||_n \geq \kappa |u|_2$, we deduce that
\begingroup
\allowdisplaybreaks
\begin{align*}
    ||\Xb u||_n^2
    % &\leq \frac{1}{n}
    % \left \langle \Xb^T \varepsilon \; , \; u \right \rangle
    % + \frac{\lambda}{n}
    % \sqrt{n} |Y - \Xb \hat \beta|_2 |u|_1 \\
    &\leq (4 + \sqrt{2}) \sigma \sqrt{\frac{\log(2p/s)}{n}}
    \left( \sqrt{s}|u|_2 + \sum_{j=s+1}^p |u|_{(j)} \right)
    + \left( \frac{|\varepsilon|_2}{\sqrt{n}}
    + ||\Xb u||_n \right)
    \frac{8}{3} \lambda \sqrt{s} |u|_2 \\
    % &\leq (4 + \sqrt{2})\sigma \sqrt{\frac{\log(2p/s)}{n}} \left( \sqrt{s}|u|_2 + |u|_1 \right)
    % + \left( 2 \sigma
    % + ||\Xb u||_n \right)
    % \frac{8}{3} \lambda \sqrt{s} |u|_2 \\
    % &\leq (4 + \sqrt{2})\sigma \sqrt{\frac{\log(2p/s)}{n}} \left( \sqrt{s}|u|_2 + \frac{8}{3} \sqrt{s} |u|_2 \right)
    % + \left( 2 \sigma
    % + ||\Xb u||_n \right)
    % \frac{8}{3} \lambda \sqrt{s} |u|_2 \\
    % &\leq (4 + \sqrt{2})\frac{11}{3}\sigma \sqrt{s\frac{\log(2p/s)}{n}} |u|_2 
    % + \left( 2 \sigma + ||\Xb u||_n \right)
    % \frac{8}{3} \lambda \sqrt{s} |u|_2 \\
    &\leq (4 + \sqrt{2})\frac{11}{3}\sigma \sqrt{s\frac{\log(2p/s)}{n}}
    \frac{||\Xb u||_n}{\kappa}
    + \left( 2 \sigma + ||\Xb u||_n \right)
    \frac{8}{3} \lambda \sqrt{s} \frac{||\Xb u||_n}{\kappa}.
\end{align*}
\endgroup
Thus,
\begin{align*}
    ||\Xb u||_n
    &\leq (4 + \sqrt{2})\frac{11}{3} \sigma
    \sqrt{s \frac{\log(2p/s)}{n}} \frac{1}{\kappa}
    + \left( 2 \sigma + ||\Xb u||_n \right)
    \frac{8}{3} \lambda \sqrt{s} \frac{1}{\kappa}.
\end{align*}
% and we get
% \begin{align*}
%     \left(1 - \frac{8 \lambda \sqrt{s}}{3 \kappa} \right)
%     ||\Xb u||_n
%     &\leq (4 + \sqrt{2})\frac{11}{3 \kappa} \sigma
%     \sqrt{s\frac{\log(2p/s)}{n}}
%     + \frac{16 \sigma \lambda \sqrt{s}}{3 \kappa}.
% \end{align*}
%
Under assumptions (\ref{cond:lambda_sqrtLasso}) and (\ref{cond:sqrt_Lasso_2}), we have
\begin{align*}
    \frac{8 \lambda \sqrt{s}}{3 \kappa}
    &=  \frac{8 \gamma}{3 \kappa}
    \sqrt{\frac{s}{n} \log \left( \frac{2p}{s} \right)}
    \leq \frac{1}{2}.
\end{align*}
Thus, we have
\begin{align}
    ||\Xb u||_n
    &\leq 2 \left( \dfrac{44 + 11\sqrt{2}}{3 \kappa} \sigma
    \sqrt{\dfrac{s}{n} \log \left( \dfrac{2p}{s} \right)}
    + \dfrac{16 \sigma \lambda \sqrt{s}}{3 \kappa}
    \right) \nonumber \\
    &\leq \dfrac{88 + 22\sqrt{2} + 32 \gamma}{3 \kappa} \sigma
    \sqrt{\dfrac{s}{n} \log \left( \dfrac{2p}{s} \right)}.
    \label{ineq_sql_n_1case}
\end{align}

We have proved in (\ref{condition_for_SRE}) that
$|u|_1 \leq (1+c_0) \sqrt{s} |u|_2$,
with $c_0 = 5/3$, so we get that
$|u|_2 \leq ||\Xb u||_n / \kappa$.
Therefore, we can deduce the following inequalities
\begin{align}
    |u|_2 &\leq \dfrac{88 + 22\sqrt{2} + 32 \gamma}{3 \kappa^2} \sigma
    \sqrt{\dfrac{s}{n} \log \left( \dfrac{2p}{s} \right)}, 
    \label{ineq_sql_2_1case} \\
    |u|_1 &\leq \dfrac{704 + 176\sqrt{2} + 256 \gamma}{9 \kappa^2} \sigma
    s \sqrt{\dfrac{1}{n} \log \left( \dfrac{2p}{s} \right)}.
    \label{ineq_sql_1_1case}
\end{align}

\vspace{0.6cm}

\noindent
{\it Second case :} $F(u) \leq G(u).$

\medskip

\noindent
Then we have
\begin{equation*}
    (4 + \sqrt{2}) \sqrt{\frac{\log(2p/s)}{n}} \left( \sqrt{s}|u|_2 + \sum_{j=s+1}^p |u|_{(j)} \right)
    \leq
    (4 + \sqrt{2}) \sqrt{\frac{\log(1/\delta_0)}{n}} ||\Xb u||_n.
\end{equation*}
Thus
\begin{equation*}
    |u|_1 \leq
    \sqrt{s}|u|_2 + \sum_{j=s+1}^p |u|_{(j)}
    \leq
    \sqrt{\frac{\log(1/\delta_0)}{\log(2p/s)}} ||\Xb u||_n.
\end{equation*}
From Lemmas \ref{first_order_consequence_sql} and \ref{lemma_Th_4_1_BLT}, we find
\begingroup
\allowdisplaybreaks
\begin{align*}
    ||\Xb u||_n^2
    % &\leq \frac{1}{n}
    % \left \langle \Xb^T \varepsilon \; , \; u \right \rangle
    % + \lambda \sqrt{n} |Y - \Xb \hat \beta|_2
    % \frac{1}{n} |u|_1 \\
    &\leq (4 + \sqrt{2}) \sigma \sqrt{\frac{\log(1/\delta_0)}{n}} ||\Xb u||_n
    + \lambda \left( \frac{|\varepsilon|_2}{\sqrt{n}}
    + ||\Xb u||_n \right) |u|_1 \\
    % &\leq (4 + \sqrt{2}) \sigma \sqrt{\frac{\log(1/\delta_0)}{n}} ||\Xb u||_n
    % + \lambda \left( \frac{||\varepsilon||_2}{\sqrt{n}}
    % + ||\Xb u||_n \right)
    % \sqrt{\frac{\log(1/\delta_0)}{\log(2p/s)}}
    % ||\Xb u||_n \\
    &\leq (4 + \sqrt{2}) \sigma \sqrt{\frac{\log(1/\delta_0)}{n}} ||\Xb u||_n
    + \lambda \left( 2 \sigma
    + ||\Xb u||_n \right)
    \sqrt{\frac{\log(1/\delta_0)}{\log(2p/s)}}
    ||\Xb u||_n.
\end{align*}
\endgroup
Thus,
\begin{align*}
    ||\Xb u||_n
    &\leq (4 + \sqrt{2}) \sigma \sqrt{\frac{\log(1/\delta_0)}{n}}
    +  \lambda \left( 2 \sigma + ||\Xb u||_n \right)
    \sqrt{\frac{\log(1/\delta_0)}{\log(2p/s)}}.
\end{align*}
We have chosen 
$\lambda = \gamma \sqrt{\frac{1}{n} \log \left( \frac{2p}{s} \right)}$, {\color{black} therefore we have
\begin{align*}
    ||\Xb u||_n
    &\leq \sigma \sqrt{\frac{\log(1/\delta_0)}{n}}
    (4 + \sqrt{2} + 2 \gamma) 
    + ||\Xb u||_n \gamma \sqrt{\frac{\log(1/\delta_0)}{n}}.
\end{align*}
By assumption, $\exp(-n/4\gamma^2) \leq \delta_0$, thus we have
\begin{align}
    ||\Xb u||_n
    &\leq \sigma \sqrt{\frac{\log(1/\delta_0)}{n}}
    (8 + 2\sqrt{2} + 4 \gamma).
    \label{ineq_sql_n_2case}
\end{align}
As a consequence, we have
\begin{align}
    |u|_1 &\leq
    \sqrt{\frac{\log(1/\delta_0)}{\log(2p/s)}} ||\Xb u||_n
    \leq \sigma \sqrt{\frac{\log^2(1/\delta_0)}{n \log(2p/s)}}
    (8 + 2\sqrt{2} + 4 \gamma).
    \label{ineq_sql_1_2case}
\end{align}
We have also
$\sqrt{s}|u|_2 \leq \sqrt{\frac{\log(1/\delta_0)}{\log(2p/s)}} ||\Xb u||_n$,
thus
\begin{equation}
    |u|_2
    \leq \sigma \sqrt{\frac{\log^2(1/\delta_0)}{s n \log(2p/s)}}
    (8 + 2\sqrt{2} + 4 \gamma).
    \label{ineq_sql_2_2case}
\end{equation}

}

As a conclusion, we can prove the result~(\ref{ineq_sql_n}) by combining the inequalities~(\ref{ineq_sql_n_1case}) and~(\ref{ineq_sql_n_2case}).
The general bound for $|u|_q$,
with $1 \leq q \leq 2$ is a consequence of the norm interpolation inequality
$|u|_q \leq |u|_1^{2/q-1} |u|_2^{2-2/q}$
which proves (\ref{ineq_sql_q}).

\begin{flushright}
    $\Box$
\end{flushright}

{\color{black}

\subsection{Proofs of the adaptive procedure}
\label{proof:adaptive_procedure}

\subsubsection{Proof of Theorem \ref{th:proc_Lepski_adaptivity}}
\label{proof:th:proc_Lepski_adaptivity}

We choose $s \in [1, s_*]$ and assume that $\beta^* \in B_0(s)$. Define $\PP := \PP_{\beta^*}$ and $m_0 := \lfloor \log_2(s) \rfloor + 1$.

\medskip

For any $a > 0$, we have 
\begin{equation}
    \PP \big( d(\tilde \beta, \beta^*) \geq a \big)
    \leq \PP \big( d(\tilde \beta, \beta^*) \geq a, \tilde m \leq m_0 \big)
    + \PP( \tilde m \geq m_0 + 1).
    \label{eq:proba_decomp_risk_tilde_beta}
\end{equation}
On the event $\{ \tilde m \leq m_0 \}$, we have the decomposition 
\begin{equation}
    d(\tilde \beta, \beta^*)
    \leq \sum_{k=\tilde m + 1}^{m_0}
    d \left(\hat \beta_{(2^{k-1})} , \hat \beta_{(2^{k})} \right)
    + d \big( \hat \beta_{(2^{m_0})}, \beta^* \big).
    \label{eq:proba_decomp_distance_beta_tilde_beta}
\end{equation}
Using Assumption \ref{assumpt:function_w}, we get that,
\begin{align}
    \sum_{k=\tilde m + 1}^{m_0}
    d \left(\hat \beta_{(2^{k-1})} , \hat \beta_{(2^{k})} \right)
    \leq \sum_{k=\tilde m + 1}^{m_0} 4 \hat \sigma C_0 w(2^k)
    \leq 4 \hat \sigma C_0 C' w(2^{m_0})
    \leq 4 \hat \sigma C_0 C' C'' w(s).
    \label{eq:bound_sum_distance_hat_beta}
\end{align}
We have $2^{m_0} \leq 2s$, therefore applying Assumption (\ref{eq:general_sup_2s}), we have with $\PP_{\beta^*}$-probability at least
$1 - \left( 2s/p \right)^{2s} - u_n$,
\begin{align}
    d(\hat \beta_{(2^{m_0})}, \beta^*)
    \leq \frac{C_2(\tilde \gamma)}{\kappa^2} \sigma w(2s)
    \leq \frac{C_2(\tilde \gamma) C''}{\kappa^2} \sigma w(s).
    \label{eq:bound_distance_beta_hat_beta}
\end{align}
Combining equations (\ref{eq:proba_decomp_distance_beta_tilde_beta}), (\ref{eq:bound_sum_distance_hat_beta}), (\ref{eq:bound_distance_beta_hat_beta}) and Assumption \ref{assumpt:bound_sigma_hat}, we get with $\PP_{\beta^*}$-probability
at least $1 - (2s/p)^{2s} - u_n - u_{n,p,M}$,
\begin{align}
    d(\tilde \beta, \beta^*)
    \leq \left( 4 \sigma C_0 C' C'' \alpha
    + \frac{C_2(\tilde \gamma) C''}{\kappa^2} \right) \sigma w(s).
    \label{eq:bound_tilde_beta_true_beta}
\end{align}

\medskip

We now bound the probability $\PP( \tilde m \geq m_0 + 1)$.
\begin{align*}
    \PP( \tilde m \geq m_0 + 1)
    &\leq \sum_{m = m_0+1}^M \PP( \tilde m = m_0 + 1)
    \leq \sum_{m = m_0+1}^M \sum_{k=m}^M \PP \bigg( 
    d \Big( \hat \beta_{(2^{k-1})} , \hat \beta_{(2^{k})} \Big) 
    > 4 \hat \sigma C_0 w(2^k) \bigg) \\
    &\leq \sum_{m = m_0+1}^M \sum_{k=m}^M
    \PP \bigg(  d \Big( \hat \beta_{(2^{k-1})} ,
    \beta^* \Big) > 2 \hat \sigma C_0 w(2^k) \bigg) + 
    \PP \bigg(  d \Big( \hat \beta_{(2^{k})} ,
    \beta^* \Big) > 2 \hat \sigma C_0 w(2^k) \bigg) \\
    &\leq 2 \sum_{m = m_0+1}^M \sum_{k=m-1}^M
    \PP \bigg(  d \Big( \hat \beta_{(2^{k-1})} ,
    \beta^* \Big) > 2 \hat \sigma C_0 w(2^k) \bigg) \\
    &\leq 2 \sum_{m = m_0+1}^M \sum_{k=m-1}^M
    \PP \bigg(  d \Big( \hat \beta_{(2^{k-1})} ,
    \beta^* \Big) > 2 \hat \sigma C_0 w(2^k) ,
    \hat \sigma \geq \frac{\sigma}{2} \bigg)
    + \PP \bigg(\hat \sigma < \frac{\sigma}{2} \bigg).
\end{align*}
Combining the previous equation with Assumption \ref{assumpt:bound_sigma_hat}, and then with Assumption (\ref{eq:general_sup_2s}), we get
\begingroup
\allowdisplaybreaks
\begin{align*}
    \PP( \tilde m \geq m_0 + 1)
    &\leq 2 \sum_{m = m_0+1}^M \sum_{k=m-1}^M
    \PP \bigg(  d \Big( \hat \beta_{(2^{k-1})} , \beta^* \Big) >
    \sigma C_0 w(2^k) \bigg) - u_{n,p,M} \\
    &\leq 2 M^2 \left( \left( \frac{2s}{p} \right)^{2s} + u_n \right) - u_{n,p,M} \\
    &\leq 2 (\log_2(s_*) + 1)^2 \left( \left( \frac{2s}{p} \right)^{2s}
    + u_n \right) - u_{n,p,M}.
\end{align*}
\endgroup
As a consequence, we deduce the bound on $\tilde s$.
Combining the last equation with equations (\ref{eq:proba_decomp_risk_tilde_beta}) and (\ref{eq:bound_tilde_beta_true_beta}), we finally get that
\begin{align*}
    &\PP \left( d(\tilde \beta, \beta^*)
    \geq \left( 4 \sigma C_0 C' C'' \alpha
    + \frac{C_2(\tilde \gamma) C''}{\kappa^2} \right) \sigma w(s) \right) \\
    & \hspace{5cm} \leq 3 (\log_2(s_*) + 1)^2 \left( \left( \frac{2s}{p} \right)^{2s}
    + u_n \right) - 2 u_{n,p,M}.
\end{align*}

\begin{flushright}
    $\Box$
\end{flushright}

\subsubsection{Proof of Lemma \ref{lemma:choice_function_w}}
\label{proof:lemma:choice_function_w}

Now, we consider the general case of the function $w(b) = b^{1/q} \sqrt{(1/n) \log(ap/b)}$, with $q$ a fixed number of the interval $[1, 2]$. The first case will correspond to $a=1$ and $q=2$ and the second case will correspond to $a=2$ with any choice of $q$.

We want to that the first part of Assumption \ref{assumpt:function_w} is satisfied, i.e., $w$ is increasing on the interval $[1,s_*]$.
Let $b \in [1,s_*]$. We have
\begin{align*}
    w'(b)
    &= \frac{1}{q} b^{(1/q)-1} \sqrt{\frac{1}{n} \log \left( \frac{ap}{b} \right) } 
    + b^{(1/q)} \frac{ -\frac{1} {nb} }
    {2 \sqrt{ \frac{1}{n} \log \left( \frac{ap}{b} \right) } } \\
    &= \frac{b^{(1/q)-1} n^{-1/2} \left( (2/q) \log \left( \frac{ap}{b} \right) - 1\right)}
    {2 \sqrt{\log \left( \frac{ap}{b} \right) } },
\end{align*}
which is positive when $(2/q) \log \left( \frac{ap}{b} \right) - 1 \geq 0$,
that is, when $b \leq ape^{-q/2}$.

We have $b \leq s_* \leq p/e = ape^{-q/2}$ when $a=1$ and $q=2$.
When $a=2$ and $q \in [1,2]$, $p/e \leq 2pe^{-1} \leq ape^{-q/2}$.
In the two cases we consider, we have proved that $w'(\cdot) \geq 0$ on the interval $[1, s_*]$, thus the function $w$ is increasing on this interval.
This proves that the first part of Assumption \ref{assumpt:function_w} is satisfied.

\medskip

Let $m$ be an integer in the interval $[1, M]$.
\begingroup \allowdisplaybreaks
\begin{align*}
    \sum_{k=1}^m w(2^k)
    &= \sum_{k=1}^m 2^{k/q} \sqrt{\frac{1}{n} \log \left(\frac{ap}{2^k} \right)} 
    = \sum_{k=0}^{m-1} 2^{(m-k)/q} \sqrt{\frac{1}{n} \log \left(\frac{ap}{2^{m-k}} \right)}  \\
    &= \frac{2^{m/q}}{\sqrt{n}} \sum_{k=0}^{m-1} \frac{1}{2^{k/q}} \sqrt{
    \left( \log \left(\frac{ap}{2^{m}} \right) + k \log(2) \right) } \\
    &\leq \frac{2^{m/q}}{\sqrt{n}} \bigg( \sum_{k=0}^{m-1} \frac{1}{2^{k/q}}
    \sqrt{\log \left(\frac{ap}{2^{m}}\right) }
    + \sum_{k=0}^{m-1} \frac{\sqrt{k}}{2^{k/q}} \sqrt{\log(2)} \bigg) \\
    &\leq \frac{2^{m/q}}{\sqrt{n}}
    \bigg(\sqrt{\log \left(\frac{ap}{2^{m}} \right)} \frac{1}{1-2^{-1/q}}
    + \sum_{k=0}^{m-1} \frac{4}{2^{k/2q}} \sqrt{\log(2)} \bigg) \\
    &\leq 2^{m/q} \sqrt{\frac{1}{n} \log \left(\frac{ap}{2^{m}} \right)} 
    \bigg( \frac{1}{1-2^{-1/q}} + \frac{4 \sqrt{\log(2)}}{1-2^{-1/(2q)}} \bigg),
\end{align*}
\endgroup
which proves that the second part is satisfied.

\medskip

\noindent
Let $b$ be an integer of $[1, s_*]$. We have
$w(2b) = (2b)^{1/q} \sqrt{(1/n) \log(2p/(2b))}
%\leq 2^{1/q} b^{1/q} \sqrt{(b/n) \log(2p/b)}
\leq 2^{1/q} w(b)$, which proves that the third part is satisfied.

\begin{flushright}
    $\Box$
\end{flushright}

\subsubsection{Proof of Lemma \ref{lemma:validity_assumption_hat_sigma}}
\label{proof:lemma:validity_assumption_hat_sigma}

We have $\beta^* \in B_0(s) \subset B_0(2^{M+1})$, therefore, we can apply Corollary \ref{cor:minimax_optimality_SQL} and Lemma \ref{lemma_epsilon_sigma}, we have with $\PP_{\beta^*}$-probability at least $1 - (2^{M+1}/p)^{2^{M+1}} - (1+e^2)e^{-n/24}$
\begin{align*}
    \hat \sigma
    &\leq ||\varepsilon||_n +
    \big|\big|\Xb(\hat \beta_{(2^{M+1})} - \beta^*) \big|\big|_n
    \leq 2 \sigma + \frac{C_2(\gamma)}{\kappa_*^2} \sigma 
    \sqrt{\frac{2^{M+1}}{n} \log \left(\frac{p}{2^{M+1}} \right)} \\
    &\leq \sigma \left(2 + \frac{C_2(\gamma)}{\kappa_*^2}
    \sqrt{\frac{2s}{n} \log \left(\frac{2p}{s} \right)} \right)
    \leq \sigma \left(2 + \frac{3\sqrt{2} C_2(\gamma)}{16 \kappa_* \gamma} \right),
\end{align*}
\begin{align*}
    \hat \sigma
    &\geq ||\varepsilon||_n -
    \big|\big|\Xb(\hat \beta_{(2^{M+1})} - \beta^*) \big|\big|_n
    \geq \frac{\sigma}{\sqrt{2}} - \frac{C_2(\gamma)}{\kappa_*^2} \sigma 
    \sqrt{\frac{2^{M+1}}{n} \log \left(\frac{p}{2^{M+1}} \right)} \\
    &\geq \sigma \left(\frac{1}{\sqrt{2}} - \frac{\sqrt{2} C_2(\gamma)}{\kappa_*^2} \sqrt{\frac{s}{n} \log \left(\frac{2p}{s} \right)} \right) \\
    &\geq \sigma \left(\frac{1}{\sqrt{2}} - \frac{\sqrt{2} C_2(\gamma)}{\kappa_*^2}
    \sqrt{\frac{2s_*}{n} \log \left(\frac{p}{s_*} \right)} \right) \\
    &\geq \sigma \left(\frac{1}{\sqrt{2}}
    - \sqrt{\left( \frac{1}{\sqrt{2}} - \frac{1}{2} \right)^2} \right)
    \geq \frac{\sigma}{2}.
\end{align*}

\begin{flushright}
    $\Box$
\end{flushright}

}

\subsection{Proof of Theorem \ref{th:minimaxSqrtSlope} }
\label{proof:minimaxSqrtSlope}

We act as in Section \ref{proof:minimaxSqrtLasso}, with suitable modifications.
We place ourselves in the event where both Lemmas \ref{lemma_epsilon_sigma} and \ref{lemma_Th_4_1_BLT} are valid, and set now $u := \hat \beta^{SQS} - \beta^*$.
Applying Lemma \ref{lemma_Th_4_1_BLT}, we will distinguish between the two cases : $G(u) \leq H(u) + \sigma |u|_2 \sqrt{\sum_{j=1}^s \lambda_j^2}$ and $H(u) + \sigma |u|_2 \sqrt{\sum_{j=1}^s \lambda_j^2} < G(u)$.

\bigskip

\noindent
{\it First case :} $G(u) \leq H(u) + \sigma |u|_2 \sqrt{\sum_{j=1}^s \lambda_j^2}.$

\medskip

\noindent
% From Lemma \ref{equ_beta_sqs}, we have
% \begin{align*}
%     \sum_{j=s+1}^p \lambda_j |u|_{(j)}
%     & \leq \sum_{j=1}^s \lambda_j |u|_{(j)}
%     + \frac{1}{\sqrt{n} | \varepsilon |_2}
%     \left \langle \Xb^T \varepsilon \; ; \; \hat \beta^{SQL} - \beta^* \right \rangle.
% \end{align*}
% %
Applying Lemma \ref{equ_beta_sqs}, Lemma \ref{lemma_Th_4_1_BLT} and then Lemma \ref{lemma_epsilon_sigma}, we have 
\begingroup \allowdisplaybreaks
\begin{align*}
    |u|_* = \sum_{j=1}^p \lambda_j |u|_{(j)}
    & \leq 2 \sum_{j=1}^s \lambda_j |u|_{(j)}
    + \frac{1}{\sqrt{n} | \varepsilon |_2}
    \left \langle \Xb^T \varepsilon \; , \; \hat \beta^{SQS} - \beta^* \right \rangle \\
    & \leq 2 \sqrt{\sum_{j=1}^s \lambda_j^2} |u|_2
    + \frac{n}{\sqrt{n} | \varepsilon |_2}
    \left((4 + \sqrt{2}) \frac{\sigma}{\gamma'} |u|_*
    + \sigma |u|_2 \sqrt{\sum_{j=1}^s \lambda_j^2} \right) \\
    & \leq 4 \sqrt{\sum_{j=1}^s \lambda_j^2} |u|_2
    + \frac{8 + 2\sqrt{2}}{\gamma'} |u|_*,
\end{align*}
\endgroup
and we get
\begin{align*}
    |u|_* \leq \dfrac{4 |u|_2}
    { 1 - \dfrac{8 + 2\sqrt{2}}{\gamma'} }
    \sqrt{\sum_{j=1}^s \lambda_j^2},
\end{align*}
Using assumption (\ref{cond:sqrt_Slope_2}), we have $\gamma' \geq 16 + 4\sqrt{2}$, therefore $|u|_* \leq 8 |u|_2 \sqrt{\sum_{j=1}^s \lambda_j^2}$. As a consequence, we get $u \in C_{WRE}(s,c_0)$ with $c_0:=8$.
Invoking Lemmas \ref{first_order_consequence_sqs}, \ref{lemma:bound_sum_lambda}, \ref{lemma_Th_4_1_BLT} and using the $WRE(s,c_0)$ condition, we get
\begingroup \allowdisplaybreaks
\begin{align*}
    ||\Xb u||_n^2
    & \leq \frac{1}{n} \left \langle \Xb^T \varepsilon \; , \; u \right \rangle
    + \frac{1}{\sqrt{n}}  |Y - \Xb \hat \beta|_2 |u|_* \\
    %
    % & \leq (4 + \sqrt{2}) \frac{\sigma}{\gamma'} |u|_*
    % + \sigma |u|_2 \sqrt{\sum_{j=1}^s \lambda_j^2}
    % + (2 \sigma + ||\Xb u||_n) |u|_* \\
    %
    & \leq (4 + \sqrt{2}) \frac{\sigma}{\gamma'} |u|_*
    + \sigma |u|_2 \sqrt{\sum_{j=1}^s \lambda_j^2}
    + (2 \sigma + ||\Xb u||_n) |u|_*  \\
    %
    % & \leq \left( (4 + \sqrt{2}) \frac{\sigma}{\gamma'} 
    % + 2 \sigma + ||\Xb u||_n \right) |u|_* 
    % + \sigma |u|_2 \sqrt{\sum_{j=1}^s \lambda_j^2} \\
    %
    & \leq \left( (32 + 8\sqrt{2}) \frac{\sigma}{\gamma'} 
    + 17 \sigma + 8 ||\Xb u||_n \right)
    |u|_2 \sqrt{\sum_{j=1}^s \lambda_j^2} \\
    %
    % & \leq \left( (32 + 8\sqrt{2}) \frac{\sigma}{\gamma'} 
    % + 17 \sigma + 8 ||\Xb u||_n \right)
    % |u|_2 \sqrt{(s/n)\log(2ep/s)} \gamma' \\
    %
    & \leq \left( (32 + 8\sqrt{2}) \frac{\sigma}{\gamma'} 
    + 17 \sigma + 8 ||\Xb u||_n \right)
    \frac{||\Xb u||_n}{\kappa'} \gamma' \sqrt{(s/n)\log(2ep/s)} .
\end{align*}
\endgroup  
Thus,
\begin{align*}
    ||\Xb u||_n
    &\leq \frac{\sigma}{\kappa'}
    \sqrt { \dfrac{s}{n} \log \left( \dfrac{2ep}{s} \right) } 
    \dfrac{32 + 8\sqrt{2} + 17 \gamma'}
    {1 - \dfrac{8 \gamma'}{\kappa'}
    \sqrt{\dfrac{s}{n}\log \left(\dfrac{2ep}{s} \right)}}.
\end{align*}
Applying condition (\ref{cond:sqrt_Slope_2}), we obtain
\begin{align}
    ||\Xb u||_n
    & \leq (64 + 16\sqrt{2} + 34 \gamma') \frac{\sigma}{\kappa'}
    \sqrt { \dfrac{s}{n} \log \left( \dfrac{2ep}{s} \right) }.
\end{align}
This and the $WRE$ condition imply
\begin{align}
    |u|_2
    & \leq (64 + 16\sqrt{2} + 34 \gamma') \frac{\sigma}{\kappa'{}^2}
    \sqrt { \dfrac{s}{n} \log \left( \dfrac{2ep}{s} \right) }.
\end{align}
Therefore, using the inequality
$|u|_* \leq 8 |u|_2
\sqrt{\sum_{j=1}^s \lambda_j^2}$, we get from Lemma \ref{lemma:bound_sum_lambda}
\begin{align}
    |u|_*
    & \leq 8(64 + 16\sqrt{2} + 34 \gamma')\gamma' \frac{\sigma}{\kappa'{}^2}
    \dfrac{s}{n} \log \left( \dfrac{2ep}{s} \right).
\end{align}

\vspace{0.6cm}

\noindent
{\it Second case :} $H(u) + \sigma |u|_2 \sqrt{\sum_{j=1}^s \lambda_j^2} \leq G(u)$.

\medskip

\noindent
Then we have
\begin{equation*}
    (4 + \sqrt{2}) \frac{\sigma}{\gamma'} |u|_* + \sigma |u|_2 \sqrt{\sum_{j=1}^s \lambda_j^2}
    \leq
    (4 + \sqrt{2}) \sigma
    \sqrt{\frac{\log(1/\delta_0)}{n}} ||\Xb u||_n.
\end{equation*}
Therefore we have
\begin{equation}
    |u|_* \leq
    \gamma' \sqrt{\frac{\log(1/\delta_0)}{n}} ||\Xb u||_n,
    \quad \text{ and } \quad
    |u|_2 \sqrt{\sum_{j=1}^s \lambda_j^2} \leq
    (4 + \sqrt{2}) \sqrt{\frac{\log(1/\delta_0)}{n}} ||\Xb u||_n.
    \label{eq:csq_2nd_case}
\end{equation}
Invoking Lemmas \ref{first_order_consequence_sqs} and \ref{lemma_Th_4_1_BLT}, and using (\ref{eq:csq_2nd_case}), we get
% \begin{align*}
%     |\Xb u|_2^2
%     % & \leq \left \langle \Xb^T \varepsilon \; , \; u \right \rangle
%     % + \sqrt{n} |Y - \Xb \hat \beta|_2 |u|_* \\
%     %
%     & \leq n (4 + \sqrt{2}) \sigma
%     \sqrt{\frac{\log(1/\delta_0)}{n}} ||\Xb u||_n
%     + n  \sigma |u|_2 \sqrt{\sum_{j=1}^s \lambda_j^2}
%     + n (2 \sigma + ||\Xb u||_n) |u|_*.
% \end{align*}
% Therefore, we have
\begin{align*}
    ||& \Xb u||_n^2
    \leq (4 + \sqrt{2}) \sigma
    \sqrt{\frac{\log(1/\delta_0)}{n}} ||\Xb u||_n
    +  \sigma |u|_2 \sqrt{\sum_{j=1}^s \lambda_j^2}
    + (2 \sigma + ||\Xb u||_n) |u|_* \\
    & \leq (4 + \sqrt{2}) \sigma
    \sqrt{\frac{\log(1/\delta_0)}{n}} ||\Xb u||_n
    + \sigma (4 + \sqrt{2}) \sqrt{\frac{\log(1/\delta_0)}{n}} ||\Xb u||_n
    + (2 \sigma + ||\Xb u||_n) \gamma'
    \sqrt{\frac{\log(1/\delta_0)}{n}} ||\Xb u||_n.
\end{align*}
which yields
\begin{align*}
    ||\Xb u||_n
    & \leq (8 + 2\sqrt{2} + 2 \gamma') \sigma
    \sqrt{\frac{\log(1/\delta_0)}{n}} 
    + ||\Xb u||_n
    \gamma' \sqrt{\frac{\log(1/\delta_0)}{n}},
\end{align*}
{\color{black}We have chosen $\exp(-n/4\gamma'{}^2) \leq \delta_0$,
which implies that
\begin{align}
    ||\Xb u||_n
    & \leq (16 + 4\sqrt{2} + 4 \gamma') \sigma
    \sqrt{\frac{\log(1/\delta_0)}{n}}.
\end{align}
We can deduce from (\ref{eq:csq_2nd_case}) that
\begin{equation}
    |u|_* \leq (16 + 4\sqrt{2} + 4 \gamma') \sigma
    \gamma' \frac{\log(1/\delta_0)}{n},
\end{equation}
and combining the second part of (\ref{eq:csq_2nd_case}) with Lemma \ref{lemma:bound_sum_lambda}, we get
\begin{equation*}
    |u|_2 \gamma' \sqrt{\frac{s}{n} \log \left(\frac{p}{s} \right)}
    \leq (4 + \sqrt{2}) \sqrt{\frac{\log(1/\delta_0)}{n}} ||\Xb u||_n
    \leq (4 + \sqrt{2}) (16 + 4\sqrt{2} + 4 \gamma') \sigma
    \frac{\log(1/\delta_0)}{n}.
\end{equation*}
Finally, we get that
\begin{equation}
    |u|_2 
    \leq \frac{(4 + \sqrt{2}) (16 + 4\sqrt{2} + 4 \gamma')}{\gamma'} \sigma
    \sqrt{\frac{\log^2(1/\delta_0)}{sn\log(p/s)}}.
\end{equation}

}

\begin{flushright}
    $\Box$
\end{flushright}

\noindent
\textbf{Acknowledgement.} This work is supported by the Labex Ecodec under the grant ANR-11-LABEX-0047 from the French Agence Nationale de la Recherche. The author thanks Professor Alexandre Tsybakov for helpful comments and discussions.

%======================
%======================

% \nocite{*} % mettre tout le contenu du .bib dans la bibliographie
\bibliographystyle{abbrv}
\bibliography{minimax_sqrt_Lasso}

\end{document}